\documentclass[twoside,12pt,reqno]{amsart}

\usepackage{amsmath,amsthm,amssymb,amstext,amsfonts,amscd}
\usepackage{graphicx}
\usepackage{multirow}
\usepackage{url}
 \newtheorem{theorem}{Theorem}[section]

 \newtheorem{corollary}{Corollary}[section]

\newtheorem{remark}{Remark}[section]

\numberwithin{equation}{section}

\begin{document}

\title[Analysis of Bell Based Euler polynomial and their Application]
{Analysis of Bell Based Euler polynomial and their Application}

\author[{\bf N. U. Khan and S. Husain}]{\bf Nabiullah Khan and Saddam Husain}

\address{Nabiullah Khan: Department of Applied
Mathematics, Faculty of Engineering and Technology,
      Aligarh Muslim University, Aligarh 202002, India}
 \email{nukhanmath@gmail.com}

\bigskip
\address{Saddam Husain: Department of Applied
    Mathematics, Faculty of Engineering and Technology,
    Aligarh Muslim University, Aligarh 202002, India}
 \email{saddamhusainamu26@gmail.com}



\begin{abstract}
In the present article, we study Bell based Euler polynomial of order $\alpha$ and investigate some useful correlation formula, summation formula and derivative formula . Also, we introduce some relation of string number of the second kind. Moreover, we drive several important formulae of bell based Euler polynomial by using umbral calculus.
\end{abstract}

\maketitle

\section{\bf{Introduction}}
The polynomials and number play an important role in the multifarious area of science such as mathematics, applied science, physics and engineering sciences and some related research area involving number theory, quantum mechanics, differential equation and mathematical physics (see \cite{Benbernou,Boas,Carlitz}). Some of the most important polynomials are Bell polynomial, Genocchi polynomial, Euler polynomial, Bernaulii polynomial and Hermite polynomial.

We know that, properties of polynomials with the help of umbral calculus studied by many authors. Kim et al. \cite{kim-kim-Rim} study some properties of Euler polynomials by using umbral Calculus and give some useful identities of Euler polynomials. Kim et al. \cite{Kim-lee} study Bernoulli, Euler and Abel polynomials by using umbral calculus and gives interesting identities of it. Kim et al. \cite{Kim-Kwonb} investigate useful identities for partially degenerate Bell Number and polynomials associated with umbral calculus and Kim et al. \cite{KIm-jang} drived some important identities and properties of umbral calculus associated with degenerate orderd Bell number and polynomials.

Recently, S. Araci et al. \cite{Araci-Acikgoz} defined a Bell based Bernaoulli polynomials. Motivated by above mention work in the present article we introduce Bell based Euler polynomials and investigate some useful correlation formula, summation formula and derivative formula of Bell based Euler polynomial of order $\alpha$. Also, we acquire some implicit summation formula and some special cases of Bell based Euler polynomials. Further, we drived some interesting result of Bell based Euler polynomials associate with umbral calculus.

\section{\bf Preliminaries}
In the present paper, we take the symbols $\mathbb{Z}, \mathbb{N}, \mathbb{N}_{0}, \mathbb{R},$ and $\mathbb{C}$ to be set of integer, set of natural number, set of non negative integers, set of real number and set of complex number respectively.

The bivariate Bell polynomials are describe by the following generating function, which is defined as follows:
\begin{equation}\label{2.1}
	\sum_{n\geq0}{\mathcal{B}_n}(x;y)\frac{{t}^n}{n!}= e^{x{t}}\,e^{y(e^{{t}}-1)}.
\end{equation}

When we take $x=0$, ${\mathcal{B}_n}(0;y)={\mathcal{B}_n}(y)$ are called classical Bell polynomials (or exponential polynomials) which is describe by following generating function (see \cite{Bell,Boas,Carlitz,Kim-Kim,kim-kim-kim-kwon}) which is defined as follows:
\begin{equation}\label{2.2}
	\sum_{n\geq0}{\mathcal{B}_n}(y)\frac{{t}^n}{n!}= e^{y(e^{{t}}-1)}.
\end{equation}

If we take $y=1$ in \eqref{2.2} i.e. ${\mathcal{B}_n}(0;1)={\mathcal{B}_n}(1)= {\mathcal{B}_n}$ are called Bell number which is defined as follows (see \cite{Bell,Boas,Carlitz,Kim-Kim,kim-kim-kim-kwon}):
\begin{equation}\label{2.3}
	\sum_{n\geq0}{\mathcal{B}_n}\frac{{t}^n}{n!}= e^{(e^{{t}}-1)}.
\end{equation}

The generating function of an Euler polynomials of order $\alpha$ (see\cite{khan-usman-wasim,khan-usman-choi,khan-choi,kim-kim-Rim,Srivastava-pinter,Srivastva-garg}) is as folows:
\begin{equation}\label{2.4}
\sum_{n\geq0}\mathcal{E}_{n}^{(\alpha)}(x)\frac{{t}^n}{n!}=e^{x{t}}\left(\frac{2}{e^{{t}}+1}\right)^{\alpha}\,\,\,\,\,\,\, ( |{t}|<2\pi ).
\end{equation}

If we take $x=0$ in \eqref{2.4} i.e. $\mathcal{E}_{n}^{(\alpha)}(0)=\mathcal{E}_{n}^{(\alpha)}$ are called Euler number  (see\cite{khan-usman-wasim,khan-usman-choi,khan-choi,kim-kim-Rim,Srivastava-pinter,Srivastva-garg}) which is shown as follows:
\begin{equation}\label{2.5}
\sum_{n\geq0}\mathcal{E}_{n}^{(\alpha)}\frac{{t}^n}{n!}=\left(\frac{2}{e^{{t}}+1}\right)^{\alpha}.
\end{equation}

The generating function of second kind string polynomials $\mathcal{S}_2(n,k;x)$ and string number $\mathcal{S}_2(n,k)$ are defined as (see\cite{Bell,Boas}):
\begin{equation}\label{2.6}
\sum_{n\geq0}\mathcal{S}_{2}(n,k;x)\frac{{t}^n}{n!}=\left(\frac{e^{{t}}-1}{k!}\right)^{k}e^{{t}x}.
\end{equation}

When $x=0$ in \eqref{2.6} i.e. $\mathcal{S}_2(n,k;0)=\mathcal{S}_2(n,k)$ are called string number and defined by following exponential generating function (see\cite{Bell,Boas}):
\begin{equation}\label{2.7}
	\sum_{n\geq0}\mathcal{S}_{2}(n,k)\frac{{t}^n}{n!}=\left(\frac{e^{{t}}-1}{k!}\right)^{k}.
\end{equation}

\section{\bf Bell based Euler polynomials and Number}
In this segment, we introduced Bell based Euler polynomials of order $\alpha$ and investigate numerous correlation formulae like implicit summation formulae, derivative formulae.

For any $n\in\mathbb{N}$ and $\alpha \in \mathbb{C}$, we define Bell based Euler polynomial of order $\alpha$ as:\\
\begin{equation}\label{3.1}
	\sum_{n\geq0}{_{\mathcal{B}} \mathcal{E}_{n}^{(\alpha)}}(x;y)\frac{{t}^n}{n!}=\left(\frac{2}{e^{{t}}+1}\right)^{\alpha}e^{{x{t}}+y(e^t -1)}\,\,\,\,\,\,\, ( |{t}|<2\pi ).
\end{equation}

If x=0 and y=1 in \eqref{3.1} then we get a Bell Based Euler number of order $\alpha$, which is defined  as follows:
\begin{equation}\label{3.2}
	\sum_{n\geq0}{_{\mathcal{B}} \mathcal{E}_{n}^{(\alpha)}}\frac{{t}^n}{n!}=\left(\frac{2}{e^{{t}}+1}\right)^{\alpha}e^{(e^{t} -1)}\,\,\,\,\,\,\, ( |{t}|<2\pi ).
\end{equation}
\subsection{Special Cases:}
In this section, we introduce some special types of Bell based Euler polynomials of order $\alpha$, which is obtain by putting particular value in \eqref{3.1} and defined as follows:

\begin{enumerate}
	\item If we choose x=0 in \eqref{3.1}, we get Bell based Euler polynomials of order $\alpha$, which are an extension of Euler polynomials of order $\alpha$ defined in \eqref{2.4} as follows:
	\begin{equation*}
   	\sum_{n\geq0}{_{\mathcal{B}} \mathcal{E}_{n}^{(\alpha)}}(y)\frac{{t}^n}{n!}=\left(\frac{2}{e^{{t}}+1}\right)^{\alpha}e^{y(e^{t} -1)}.
    \end{equation*}
  
    \item In case y=0 in \eqref{3.1} the Bell based Euler polynomials of order $\alpha$ reduced to the famililar Euler polynomials $\mathcal{E}_{n}^{(\alpha)}(x)$ of order $\alpha$ defined in \eqref{2.4}
	\begin{equation*}
	\sum_{n\geq0}{_{\mathcal{B}} \mathcal{E}_{n}^{(\alpha)}}(x)\frac{{t}^n}{n!}=\left(\frac{2}{e^{{t}}+1}\right)^{\alpha}e^{x{t}}.
	\end{equation*}

    \item In case y=0 and $\alpha$=1 in \eqref{3.1} the Bell Based Euler polynomials ${_{\mathcal{B}} \mathcal{E}_{n}^{(\alpha)}}(x;y)$ reduced to usual Euler polynomials $\mathcal{E}_{n}(x)$ defined as:
    \begin{equation*}
    	\sum_{n\geq0}{_{\mathcal{B}} \mathcal{E}_{n}}(x)\frac{{t}^n}{n!}=\left(\frac{2}{e^{{t}}+1}\right)e^{x{t}}.
    \end{equation*}
 \end{enumerate}

 \begin{theorem}
 The  following relation hold true for $\alpha \in \mathbb{C}$ and $n \in \mathbb{N}$;
 
 \begin{equation}\label{3.3}
 	{_{\mathcal{B}} \mathcal{E}_{n}^{(\alpha)}}(x;y)=\sum\limits_{k=0}^{n}\binom{n}{k} { {\mathcal{E}}_{k}^{(\alpha)}}(x) {\mathcal{B}}_{n-k}(y).
 \end{equation}

\begin{proof} By using relation \eqref{3.1}, we have
	
\begin{equation*}
\aligned
\sum_{n\geq0}{_{\mathcal{B}} \mathcal{E}_{n}^{(\alpha)}}(x;y)\frac{{t}^n}{n!}=&\left(\frac{2}{e^{{t}}+1}\right)^{\alpha}e^{{x{t}}+y(e^{t} -1)}\\
=&\left\{ \left(\frac{2}{e^{{t}}+1}\right)^{\alpha}e^{x{t}}\right\} \left\{e^{y(e^{t} -1)}\right\}\\
=& \left\{\sum_{k\geq0}{\mathcal{E}_{k}^{(\alpha)}}(x)\frac{{t}^k}{k!}\right\}\left\{\sum_{n\geq0}{\mathcal{B}_n}(y)\frac{{t}^n}{n!}\right\}.
\endaligned
\end{equation*}

Now, using series rearrangement method, we get
\begin{equation*}
\aligned
\sum_{n\geq0}{_{\mathcal{B}} \mathcal{E}_{n}^{(\alpha)}}(x;y)\frac{{t}^n}{n!}=& \sum_{n\geq0}\left\{\sum\limits_{k=0}^{n}\binom{n}{k} { {\mathcal{E}}_{k}^{(\alpha)}}(x) {\mathcal{B}}_{n-k}(y)\right\}\frac{{t}^n}{n!}.
\endaligned
\end{equation*}

By equating same power of t both side, we get desired result \eqref{3.3}. 
\end{proof}
\end{theorem}

\begin{theorem}
Bell based Euler polynomials of order $\alpha$ satisfy following relation for any $\alpha \in \mathbb{C}$ and $n \in \mathbb{N}$:
\begin{equation}\label{3.4}
	{_{\mathcal{B}} \mathcal{E}_{n}^{(\alpha)}}(x;y)=\sum\limits_{k=0}^{n}\binom{n}{k} { {\mathcal{E}}_{k}^{(\alpha)}} {\mathcal{B}}_{n-k}(x;y).
\end{equation}

\begin{proof}
 By	using generating function \eqref{3.1}, we have
 \begin{equation*}
 	\aligned
 	\sum_{n\geq0}{_{\mathcal{B}} \mathcal{E}_{n}^{(\alpha)}}(x;y)\frac{{t}^n}{n!}=&\left(\frac{2}{e^{{t}}+1}\right)^{\alpha}e^{{x{t}}+y(e^{t} -1)}\\
 	=&\left\{ \left(\frac{2}{e^{{t}}+1}\right)^{\alpha}\right\} \left\{e^{x{t}+y(e^{t} -1)}\right\}\\
 	=& \left\{\sum_{k\geq0}{\mathcal{E}_{k}^{(\alpha)}}\frac{{t}^k}{k!}\right\}\left\{\sum_{n\geq0}{\mathcal{B}_n}(x;y)\frac{{t}^n}{n!}\right\}\\
 	\endaligned
 \end{equation*}

After applying series rearrangement technique and compare same power of t, we get desired result \eqref{3.4}.
\end{proof}
\end{theorem}	
	
\begin{theorem}
If $\alpha \in \mathbb{C}$ and $n \in \mathbb{N}$, following relation hold true;	
\begin{equation}\label{3.5}
	{_{\mathcal{B}} \mathcal{E}_{n}^{(\alpha)}}(x;y)=\sum\limits_{k=0}^{n}\binom{n}{k} { {\mathcal{E}}_{k}^{(\alpha)}}(y)\,\,x^{n-k}.
\end{equation}
\begin{proof}
 Using relation \eqref{3.1}, we have
 \begin{equation*}
 	\aligned
 	\sum_{n\geq0}{_{\mathcal{B}} \mathcal{E}_{n}^{(\alpha)}}(x;y)\frac{{t}^n}{n!}=&\left(\frac{2}{e^{{t}}+1}\right)^{\alpha}e^{{x{t}}+y(e^{t} -1)}\\
 	=&\left\{ \left(\frac{2}{e^{{t}}+1}\right)^{\alpha}e^{y(e^{{t}}-1)}\right\} \left\{e^{x{t}}\right\}\\
 	=& \left\{\sum_{k\geq0}{_\mathcal{B}}{\mathcal{E}_{k}^{(\alpha)}}(y)\frac{{t}^k}{k!}\right\}\left\{\sum_{n\geq0}\frac{(x{t})^n}{n!}\right\}\\
 	=& \left\{\sum_{n\geq0} \sum_{k\geq0} {_\mathcal{B}}{\mathcal{E}_{k}^{(\alpha)}}(y)\,\,\frac{x^n}{n!}\frac{{t}^{n+k}}{k!}\right\}.
 	\endaligned
 \end{equation*}
By using series rearrangement, we get 
\begin{equation*}
 	\sum_{n\geq0}{_{\mathcal{B}} \mathcal{E}_{n}^{(\alpha)}}(x;y)\frac{{t}^n}{n!}= \sum_{n\geq0}\left\{ \sum_{k=0}^{n}\binom{n}{k}{_\mathcal{B}}{\mathcal{E}_{k}^{(\alpha)}}(y)\,\,x^{n-k}\right\}\frac{{t}^{n}}{n!}.
\end{equation*}

By equating same power of t both side, we get desired result \eqref{3.5}. 
\end{proof}
\end{theorem}

\section{\bf Implicit summation formulae}
In this section, we discus some useful implicit summation formulae for Bell Based Euler polynomials of order $\alpha$, which is defined in following theorem as follows:

\begin{theorem}
If $\alpha_{1},\alpha_{2} \in \mathbb{C}$ and $n \in \mathbb{N}$, following relation hold true;
	
 \begin{equation}\label{4.1}
	{_{\mathcal{B}} \mathcal{E}_{n}^{(\alpha_{1}+\alpha_{2})}}(x_{1}+x_{2};y_{1}+y_{2})=\sum\limits_{k=0}^{n}\binom{n}{k} {_\mathcal{B} {\mathcal{E}}_{k}^{(\alpha_{1})}}(x_{1};y_{1})\,\, {_\mathcal{B}}{\mathcal{E}}^{(\alpha_{2})}_{n-k}(x_{2};y_{2}).
 \end{equation}	
	
\begin{proof}
	Using the following identity

$$\left(\frac{2}{e^{{t}}+1}\right)^{\alpha_{1}+\alpha_{2}}e^{{(x_{1}+x_{2})t}+(y_{1}+y_{2})(e^t -1)}$$

$$=\left\{ \left(\frac{2}{e^{{t}}+1}\right)^{\alpha_{1}}e^{x_{1}{t}}+y_{1}(e^{t} -1)\right\} \left\{ \left(\frac{2}{e^{{t}}+1}\right)^{\alpha_{2}}e^{x_{2}{t}}+y_{2}(e^{t} -1)\right\}.$$

 By using generating function \eqref{3.1}, we have
	
$$\sum_{n\geq0}{_{\mathcal{B}} \mathcal{E}_{n}^{(\alpha_{1}+\alpha_{2})}}(x_{1}+x_{2};y_{1}+y_{y})\frac{{t}^n}{n!}$$
\begin{equation*}
\aligned	
=&\left(\frac{2}{e^{{t}}+1}\right)^{\alpha_{1}+\alpha_{2}}e^{{(x_{1}+x_{2}){t}}+(y_{1}+y_{2})(e^{t} -1)}\\
=&\left\{ \left(\frac{2}{e^{{t}}+1}\right)^{\alpha_{1}}e^{x_{1}{t}}+y_{1}(e^{t} -1)\right\} \left\{ \left(\frac{2}{e^{{t}}+1}\right)^{\alpha_{2}}e^{x_{2}{t}}+y_{2}(e^{t} -1)\right\}\\
=&\left\{\sum_{k\geq0}{_\mathcal{B}}{\mathcal{E}_{k}^{(\alpha_{1})}}(x_{1};y_{1})\frac{{t}^k}{k!}\right\} \left\{\sum_{n\geq0}{_\mathcal{B}}{\mathcal{E}_{n}^{(\alpha_{2})}}(x_{2};y_{2})\frac{{t}^n}{n!}\right\}\\
=&\left\{\sum_{n\geq0}\sum_{k\geq0}{_\mathcal{B}}{\mathcal{E}_{k}^{(\alpha_{1})}}(x_{1};y_{1})\,\,{_\mathcal{B}}{\mathcal{E}_{n}^{(\alpha_{2})}}(x_{2};y_{2})\,\,\frac{{t}^{n+k}}{n! k!}\right\}
\endaligned	
\end{equation*}

Using series rearrangement technique, we obtain
\begin{equation*}
\sum_{n\geq0}{_{\mathcal{B}} \mathcal{E}_{n}^{(\alpha_{1}+\alpha_{2})}}(x_{1}+x_{2};y_{1}+y_{y})\frac{{t}^n}{n!}
=\sum_{n\geq0}\left\{\sum_{k=0}^{n}\binom{n}{k}{_\mathcal{B}}{\mathcal{E}_{k}^{(\alpha_{1})}}(x_{1};y_{1})\,\,{_\mathcal{B}}{\mathcal{E}_{n-k}^{(\alpha_{2})}}(x_{2};y_{2})\,\,\right\}\frac{{t}^{n}}{n!}
\end{equation*}

Now, equating same power of ${t}$ both side, we get desired result \eqref{4.1}.
\end{proof}	
\end{theorem}

\begin{remark}
In case, if we choose $\alpha_{1}=\alpha$, $\alpha_{2}=0$, $x_{1}=x$, $x_{2}=1$, $y_{1}=y$ and $y_{2}=0$ in \eqref{4.1}, we have
\begin{equation}\label{4.2}
	{_{\mathcal{B}} \mathcal{E}_{n}^{(\alpha)}}(x+1;y)=\sum\limits_{k=0}^{n}\binom{n}{k} {_\mathcal{B} {\mathcal{E}}_{k}^{(\alpha)}}(x;y),
\end{equation}
which is an extension of Euler polynomials defined by
\begin{equation}\label{4.3}
	{\mathcal{E}_{n}}(x+1)=\sum\limits_{k=0}^{n}\binom{n}{k} {\mathcal{E}_{k}}(x;y).
\end{equation} 	
\end{remark}

\begin{theorem}
If $\alpha \in \mathbb{C}$ and $n \in \mathbb{N}$, Bell bassed Euler polynomials of order $\alpha$ satisfy following summation formula:
\begin{equation}\label{4.4}
{_{\mathcal{B}} \mathcal{E}_{n+1}^{(\alpha)}}(x+1;y)-	{_{\mathcal{B}} \mathcal{E}_{n+1}^{(\alpha)}}(x;y)=\sum\limits_{k=0}^{n}\binom{n+1}{k} {_\mathcal{B} {\mathcal{E}}_{k}^{(\alpha)}}(x;y)
\end{equation}

\begin{proof}
Using relation defined in \eqref{3.1}, we get
$$\sum_{n\geq0}{_{\mathcal{B}} \mathcal{E}_{n}^{(\alpha)}}(x+1;y)\frac{{t}^n}{n!}-\sum_{n\geq0}{_{\mathcal{B}} \mathcal{E}_{n}^{(\alpha)}}(x;y)\frac{{t}^n}{n!}$$

\begin{equation*}
	\aligned
	=&\left(\frac{2}{e^{{t}}+1}\right)^{\alpha}e^{{(x+1){t}}+y(e^{t} -1)}-\left(\frac{2}{e^{{t}}+1}\right)^{\alpha}e^{{x{t}}+y(e^{t} -1)}\\
	=& \left(\frac{2}{e^{{t}}+1}\right)^{\alpha}e^{{x{t}}+y(e^{t} -1)}\,\,(e^{{t}}-1)\\
	=& \left\{\sum_{k\geq0}{_{\mathcal{B}} \mathcal{E}_{k}^{(\alpha)}}(x;y)\frac{{t}^k}{k!}\right\}\,\,\left\{\sum_{n\geq0}\frac{{t}^{n+1}}{(n+1)!}\right\}.
	\endaligned
\end{equation*}

After, series rearrangement to compare both side the same power of ${t}$, we obtain a desired result \eqref{4.4} 
\end{proof}
\end{theorem}

\begin{theorem}
If $\alpha=1$ and $n\in \mathbb{N}$ then Bell based Euler polynomials satisfy following relation:
\begin{equation}\label{4.5}
	{\mathcal{B}_n}(x;y)=\frac{{_{\mathcal{B}} \mathcal{E}_{n}}(x+1;y)+	{_{\mathcal{B}} \mathcal{E}_{n}}(x;y)}{2}.
\end{equation}

\begin{proof}
 Using generating function \eqref{3.1} for $\alpha=1$ and definition of bivariate Bell polynomials, we get
 \begin{equation*}
 	\aligned
 	\sum_{n\geq0}{\mathcal{B}_n}(x;y)\frac{{t}^n}{n!}=& e^{x{t}+{y(e^{{t}}-1)}}\\
 	=&\frac{e^{{t}}+1}{2}\,\,\left\{\sum_{k=0}^{\infty}{_{\mathcal{B}} \mathcal{E}_{k}}(x;y)\right\}\\
 	=& \frac{e^{{t}}+1}{2}\left\{\left(\frac{2}{e^{{t}}+1}\right)e^{{x{t}}+y(e^{t} -1)}\right\}\\
 	=&\frac{1}{2}\left\{\left(\frac{2}{e^{{t}}+1}\right)e^{{(x+1){t}}+y(e^{t} -1)}+\left(\frac{2}{e^{{t}}+1}\right)e^{{x{t}}+y(e^{t} -1)}\right\}\\
 	=&\frac{1}{2}\left\{\sum_{n\geq0}{_{\mathcal{B}}\mathcal{E}_{n}}(x+1;y)\frac{{t}^n}{n!}+\sum_{n\geq0}{_{\mathcal{B}} \mathcal{E}_{n}}(x;y)\frac{{t}^n}{n!}\right\}.
 	\endaligned
 \end{equation*}
By equating same power of ${t}$ both side, we get desired result \eqref{4.5}.
\end{proof}
\end{theorem}

\begin{remark}
The formula \eqref{4.5} is generalized form of Euler polynomials defined as
\begin{equation}\label{4.6}
	x^{n}=\frac{{ \mathcal{E}_{n}}(x+1)+	{ \mathcal{E}_{n}}(x)}{2}.
\end{equation}
\end{remark}

\begin{theorem}
If $n\geq 0$, then 
\begin{equation}\label{4.7}
{_{\mathcal{B}} \mathcal{E}_{n}^{(\alpha)}}(x;y)= \sum_{j=0}^{n}\sum_{k\geq0}\binom{n}{j}(x)_{k}\,\mathcal{S}_{2}(j, k)\,{_{\mathcal{B}} \mathcal{E}_{n}^{(\alpha)}}(y).
\end{equation}

\begin{proof}
 By using \eqref{3.1}, we have
 \begin{equation*}
 	\aligned
 	\sum_{n\geq0}{_{\mathcal{B}} \mathcal{E}_{n}^{(\alpha)}}(x;y)\frac{{t}^n}{n!}=&\left(\frac{2}{e^{{t}}+1}\right)^{\alpha}e^{{x{t}}+y(e^{t} -1)}\\
 	=&\left(\frac{2}{e^{{t}}+1}\right)^{\alpha}e^{y(e^{t} -1)}\,e^{x{t}}\\
 	=&\left(\frac{2}{e^{{t}}+1}\right)^{\alpha}e^{y(e^{t} -1)}\,(1+e^{{t}} -1)^x\\
 	=&\left\{\sum_{n\geq0}{_{\mathcal{B}} \mathcal{E}_{n}^{(\alpha)}}(x;y)\frac{{t}^n}{n!}\right\}\left\{\sum_{k\geq0}(x)_{k}\,\frac{(e^{{t}} -1)^{k}}{k!}\right\}\\
 	=&\left\{\sum_{n\geq0}{_{\mathcal{B}} \mathcal{E}_{n}^{(\alpha)}}(x;y)\frac{{t}^n}{n!}\right\}\left\{\sum_{k\geq0}(x)_{k}\,\sum_{j\geq0}\mathcal{S}_{2}(j, k)\frac{{t}^{j}}{j!}\right\}\\
 	\endaligned
 \end{equation*}
By using series rearrangement technique and equating same power of $\mathfrak{t}$ both side, we get desired result \eqref{4.7}.
\end{proof}
\end{theorem}

\section{\bf Deriavtive formula}
\begin{theorem}
The differential operator formula for the bell based Euler polynomials of order $\alpha$ w.r.t. x defined as follows:
\begin{equation}\label{5.1}
	\frac{\partial}{\partial x}\,\, {_{\mathcal{B}} \mathcal{E}_{n}^{(\alpha)}}(x;y)= n\,\,{_{\mathcal{B}} \mathcal{E}_{n-1}^{(\alpha)}}(x;y),
\end{equation}
hold for all $n\in\mathbb{N}$.

\begin{proof}
  we know that
  \begin{equation}\label{5.2}
  \frac{\partial}{\partial x}\,\,{e^{x{t}+y(e^{{t}}-1)}}=t\,\,e^{x{t}+y(e^{{t}}-1)}
  \end{equation}

  By using definition \eqref{3.1} in \eqref{5.2}, we obtain desired result \eqref{5.1}.
\end{proof}
\end{theorem}

\begin{theorem}
	The difference operator formula for Bell based Euler polynomials of order $\alpha$ w.r.t. y defined as:
\begin{equation}\label{5.3}
	\frac{\partial}{\partial y}\,\, {_{\mathcal{B}} \mathcal{E}_{n}^{(\alpha)}}(x;y)=(-2)\,\left\{{_{\mathcal{B}} \mathcal{E}_{n}^{(\alpha)}}(x;y)-{_{\mathcal{B}} \mathcal{E}_{n}^{(\alpha-1)}}(x;y)\right\},
\end{equation}
hold for all $n\in\mathbb{N}$.
\begin{proof}
	By using well known derivative properties
	\begin{equation}\label{5.4}
		\frac{\partial}{\partial y}\,\,{e^{x{t}+y(e^{{t}}-1)}}=(e^{{t}}-1)\,\,e^{x{t}+y(e^{{t}}-1)}
	\end{equation}

Now, using definition \eqref{3.1}, we get
 \begin{equation*}
\aligned
\frac{\partial}{\partial y}\left\{\sum_{n\geq0}{_{\mathcal{B}} \mathcal{E}_{n}^{(\alpha)}}(x;y)\frac{{t}^n}{n!}\right\}=&\frac{\partial}{\partial y}\left\{\left(\frac{2}{e^{{t}}+1}\right)^{\alpha}e^{{x{t}}+y(e^{t} -1)}\right\}\\
=&\left\{ \left(\frac{2}{e^{{t}}+1}\right)^{\alpha}e^{x{t}+y(e^{{t}}-1)}\right\}(e^{{t}}-1)\\
=&(-2+e^{{t}}+1)\,\left\{\left(\frac{2}{e^{{t}}+1}\right)^{\alpha}e^{x{t}+y(e^{{t}}-1)}\right\}\\
=&(-2)\left(1-\frac{e^{{t}}+1}{2}\right)\,\left\{\left(\frac{2}{e^{{t}}+1}\right)^{\alpha}e^{x{t}+y(e^{{t}}-1)}\right\}
\endaligned
\end{equation*}

Now, using \eqref{3.1} l.h.s of above equation and equating  same power of $\mathfrak{t}$ both side, we get desired result \eqref{5.3}.
\end{proof}	
\end{theorem}

\section{\bf Application of Bell based Euler polynomials from umbral calculus}
In this section, we discus the concept of umbral calculus (see\cite{Simsek-Srivastva,Dere-simsek,Araci-Acikgoz,kim-kim-Rim,Kim-lee,Kim-Kwonb,KIm-jang-jang}) and by using the concept of umbral calculus we drive some useful relation of Bell Based Euler polynomilas.

Let $\mathcal{G}$ be set of all formal power series in the variable t over complex field $\mathbb{C}$ with
\begin{equation}\label{6.1}
	\mathcal{G}= \left\{\mathfrak{g}|\mathfrak{g}(t)=\sum_{k\geq0}\frac{c_k}{k!}{t}^{k}\,\,s.t.\,\, c_{k}\in\mathbb{C}\right\}.
\end{equation}

Let $\mathcal{P}$ be set of polynomials in the single vriable $t$ and $\mathcal{P}^*$ be set of vector space of all linear functional on $\mathcal{P}$. In the umbral calculus, we denote $\left<\mathcal{T}|~\mathfrak{q}(x)\right>$ be linear functional $\mathcal{T}$ on the polynomials $\mathfrak{q}(x)$. Now, we defined the vector space operations on $\mathcal{P}^*$ as follows:

$$\left<\mathcal{T}_{1}+\mathcal{T}_{2}|\,\mathfrak{q}(x)\right>=\left<\mathcal{T}_{1}|\,\mathfrak{q}(x)\right>+\left<\mathcal{T}_{2}|\,\mathfrak{q}(x)\right>$$
and
$$\left<\beta \mathcal{T}|~\mathfrak{q}(x)\right>=\beta \left<\mathcal{T}|~\mathfrak{q}(x)\right>$$

for any constant $\beta$ in $\mathbb{C}$.\\

The formal power series
\begin{equation}\label{6.2}
	\mathfrak{g}(t)=\sum_{k\geq0}\frac{c_k}{k!}{t}^{k}\,\, \in \mathcal{G}
\end{equation}

defined a linear functional on $\mathcal{P}$ as
\begin{equation}\label{6.3}
	\left<\mathfrak{g}(t)|~x^n\right>=c_n
\end{equation}

for all $n \in \mathbb{N} \cup \{0\}$.\\

If we choose $\mathfrak{g}(t)={t}^{k}$ in \eqref{6.2} and \eqref{6.3}, we get
\begin{equation}\label{6.4}
	\left<{t}^{k}|~x^n\right>=n!\delta_{n,k}, 
\end{equation}

for all  $n,k \in \mathbb{N}\cup\{0\}$ and
\begin{equation*}
	\delta_{n,k} = \left\{
	\begin{array}{ll}
		0 & if\,\, n \neq k\\
		1, & if\,\, n=k
	\end{array}.
	\right.
\end{equation*}

Since any linear functional $\mathcal{T}$ in $\mathcal{P}^*$ has in the form of \eqref{6.2} i.e.,
\begin{equation*}
	\mathfrak{g}_{\mathcal{T}}({t})=\sum_{k\geq0}\left<\mathcal{T}|~x^{k}\right>\frac{{t}^{k}}{k!}
\end{equation*}

and
\begin{equation*}
	\left<\mathfrak{g}_{\mathcal{T}}({t})|~x^n\right>=\left<\mathcal{T}|~x^n\right>.
\end{equation*}

So, the linear functional $\mathcal{T}=\mathfrak{g}_{\mathcal{T}}({t})$. We know that, the map $\mathcal{T}\rightarrow \mathfrak{g}_{\mathcal{T}}({t})$ is a vector space isomorphism from set of vector space of all linear functional on $\mathcal{P}$ onto set of formal series $\mathcal{G}$. Therefore, set of formal series $\mathcal{G}$ have vector space of all linear functionals on $\mathcal{P}$,  also $\mathcal{G}$ have algebra of formal power series and so for $\mathfrak{g}(t) \in \mathcal{G}$ will be treatrd as both a formal power series  and  a linear functional. From \eqref{6.3}, we get
\begin{equation}\label{6.5}
	\left<e^{y{t}}|~x^n\right>=y^{n}
\end{equation}
and so, we have
\begin{equation}\label{6.6}
	\left<e^{y{t}}|~\mathfrak{q}(x)\right>=\mathfrak{q}(y)
\end{equation}
for all $\mathfrak{q}(x)\in \mathcal{P}$.

\vspace{0.25cm}
We know that the order of power series $\mathfrak{g}(t)$ $(i.e.~~o(\mathfrak{g}(t))$ will be the smallest positive integer $k$ such that the coefficient of ${t}^{k}$ does not vanish. We know that $\mathfrak{g}(t)$ is invertible series if the order of formal power series $\mathfrak{g}(t)$ is zero, also the formal power series $\mathfrak{g(t)}$ is a delta series if the order of $\mathfrak{g}(t)$ is one $(i.e.\, o(\mathfrak{g}(t)=1)$. (see\cite{Simsek-Srivastva,Dere-simsek,Araci-Acikgoz,kim-kim-Rim,Kim-lee,Kim-Kwonb,KIm-jang-jang}).

\vspace{0.25cm}
For $\mathfrak{g}_{1}({t}),...,\mathfrak{g}_{m}({t}) \in \mathcal{G}$, then
\begin{equation*}
	\left<\mathfrak{g}_{1}({t}),...,\mathfrak{g}_{m}({t})\arrowvert~ x^{n}\right>=\sum_{i_{1}+i_{2}+...+i_{m}=n}\binom{n}{i_{1},...,i_{m}}\left<\mathfrak{g}_{1}({t})|~x^{i_{1}}\right>...\left<\mathfrak{g}_{m}({t})|~x^{i_{m}}\right>,
\end{equation*}

where 
\begin{equation*}
	\binom{n}{i_{1},...,i_{m}}=\frac{n!}{i_{1}! ... i_{m}!}.
\end{equation*}

If $\mathfrak{g}(t), \mathfrak{h}(t) \in \mathcal{G}$, then
\begin{equation}\label{6.7}
	\left<\mathfrak{g}(t)\mathfrak{h}(t)|~\mathfrak{q}(x)\right>=\left<\mathfrak{g}(t)|~\mathfrak{h}(t)\mathfrak{q}(x)\right>=\left<\mathfrak{h}(t)|~\mathfrak{g}(t)\mathfrak{q}(x)\right>.
\end{equation}

\vspace{0.25cm}
Therfore, $\forall$ $\mathfrak{g}(t) \in \mathcal{G}$
\begin{equation}\label{6.8}
	\mathfrak{g}(t)=\sum_{k\geq0}\left<\mathfrak{g}(t)|~x^{k}\right>\frac{{t}^{k}}{k!},
\end{equation}

and $\forall$ polynomials $\mathfrak{q}(x)$
\begin{equation}\label{6.9}
	\mathfrak{q}(x)=\sum_{k\geq0}\frac{\left<{t}^k|~\mathfrak{q}(x)\right>}{k!}{t}^{k}.
\end{equation}

By using \eqref{6.9}, we get
\begin{equation}\label{6.10}
	\mathfrak{q}^{k}(x)= D^{k}\mathfrak{q}(x)=\sum_{k\geq0}\frac{\left<{t}^l|~\mathfrak{q}(x)\right>}{l!}x^{l-k}\prod_{s=1}^{k}(l-s+1)
\end{equation}

From \eqref{6.10}, we have
\begin{equation}\label{6.11}
	\mathfrak{q}^{k}(0)=\left<{t}^k|~\mathfrak{q}(x)\right>\,\,\,\, and\,\,\,\,\,\,\,\mathfrak{q}^{k}(0)=\left<1|~\mathfrak{q}^{k}(x)\right>
\end{equation}

By \eqref{6.11}, we note that
\begin{equation}\label{6.12}
	{t}^{k}\mathfrak{q}(x)=\mathfrak{q}^k(x).
\end{equation}

\vspace{0.25cm}
Let $\mathfrak{g}(t)$ and $\mathfrak{h}(t)$ be an element of formal powr series $\mathcal{G}$ such that $\mathfrak{g}(t)$ be a delta and $\mathfrak{h}(t)$ be an invertible series. Then $\exists$ an unique sequence $\mathcal{S}_{n}(x)$ of the polynomials with following properties
\begin{equation}\label{6.13}
	\left<\mathfrak{h}(t)\mathfrak{g}(t)^{k}\arrowvert ~\mathcal{S}_{n}(x)\right>=n!\delta_{n,k}
\end{equation}

For all $n,k \in \mathbb{N}\cup \{0\}$, which is orthogonality condition for the Sheffer sequence (see\cite{Simsek-Srivastva,Dere-simsek,Araci-Acikgoz,kim-kim-Rim,Kim-lee,Kim-Kwonb,KIm-jang-jang}).

\vspace{0.25cm}
The sequence $\mathcal{S}_{n}(x)$ is said to be Sheffer sequence for $(\mathfrak{h}(t), \mathfrak{g}(t))$, which is denoted by $\mathcal{S}_{n}(x)\sim (\mathfrak{h}(t), \mathfrak{g}(t))$.

Let $\mathcal{S}_{n}(x)$  is the Sheffer sequence for $(\mathfrak{h}(t), \mathfrak{g}(t))$. Then for $\mathfrak{f(t)}\in \mathcal{G}$ and for $\mathfrak{q}(x)$, we have the the following relation;
\begin{equation}\label{6.14}
	\mathfrak{f}(t)=\sum_{k\geq0}\frac{\left<\mathfrak{f}(t)\arrowvert ~\mathcal{S}_{k}(x)\right>}{k!}\mathfrak{h}(t)\mathfrak{g}(t)^{k}
\end{equation}

and
\begin{equation}\label{6.15}
	\mathfrak{q}(x)=\sum_{k\geq0}\frac{\left<\mathfrak{h}(t)\mathfrak{g}(t)^{k}\arrowvert ~\mathfrak{q}(x)\right>}{k!}\mathcal{S}_{k}(x)
\end{equation}

and the sequence $\mathcal{S}_{n}(x)$ be a Sheffer sequence for ${(\mathfrak{h}(t), \mathfrak{g}(t))}$, iff
\begin{equation}\label{6.16}
	\frac{1}{\mathfrak{h}(\bar{\mathfrak{g}}({t}))}e^{y\bar{\mathfrak{g}}({t})}=\sum_{k\geq0}\frac{\mathcal{S}_{k}(y)}{k!}{t}^{k},\,\,\,\,\,\,\,\, \forall\,\,\,y\in\mathbb{C}.
\end{equation}

Here, $\bar{\mathfrak{g}}({t})$ is composition inverse of $\mathfrak{g}(t)$ i.e. $\bar{\mathfrak{g}}(\mathfrak{g}(t))=\mathfrak{g}(\bar{\mathfrak{g}}({t}))={t}$.

\vspace{0.25cm}
Suppose that $\mathcal{S}_{n}(x)$ be an appell sequence for $\mathfrak{h}(t)$. From \eqref{6.16}, we get

\begin{equation}\label{6.17}
	\mathcal{S}_{n}(x)=\frac{1}{\mathfrak{h}(t)}x^{n} 	\Leftrightarrow t\mathcal{S}_{n}(x)=n\mathcal{S}_{n-1}(x).
\end{equation}

\vspace{0.25cm}
Recently, many author have studied Euler polynomilas, Bernoulli polynomials and Bell polynomials under theory of umbral calculus (see\cite{Simsek-Srivastva,Dere-simsek,Araci-Acikgoz,kim-kim-Rim,Kim-lee,Kim-Kwonb,KIm-jang,KIm-jang-jang}) .

Recall from \eqref{3.1}, we have
\begin{equation}\label{6.18}
	\sum_{n\geq0}{_{\mathcal{B}} \mathcal{E}_{n}^{(\alpha)}}(x;y)\frac{{t}^n}{n!}=\left(\frac{2}{e^{{t}}+1}\right)e^{{x{t}}+y(e^{t} -1)}\,\,\,\,\,\,\, ( |{t}|<2\pi ).
\end{equation}

As $t$ goes to zero in \eqref{6.18} gives ${_{\mathcal{B}} \mathcal{E}_{n}}(x;y)$ equal to one $(i.e.\,\,{_{\mathcal{B}} \mathcal{E}_{n}}(x;y)=1 )$ which means that $o\left(\left(\frac{2}{e^{{t}}+1}\right)e^{{x{t}}+y(e^{t} -1)}\right)=0$, which implies that \eqref{3.1} is an invertible series and treated as a Sheffer sequence.

We defined an important and useful properties of Bell based Euler polynomials by using the definition of umbral calculus. Which is defined as follows:

From \eqref{6.16} and \eqref{6.18}, we have
\begin{equation}\label{6.19}
	{_{\mathcal{B}} \mathcal{E}_{n}}(x;y) \sim \left(\frac{e^{{t}}+1}{2}\,e^{-y(e^{{t}}-1)}, {t}\right)
\end{equation}
and 
\begin{equation}\label{6.20}
	{t}\,\,{_{\mathcal{B}} \mathcal{E}_{n}}(x;y)=n\,\,{_{\mathcal{B}} \mathcal{E}_{n-1}}(x;y)
\end{equation}

\vspace{0.25cm}
Then from \eqref{6.19} and \eqref{6.20} we say that ${_{\mathcal{B}} \mathcal{E}_{n}}(x;y)$  is an appell sequence for $\frac{e^{{t}}+1}{2}\,e^{-y(e^{{t}}-1)}$.

\begin{theorem}
 If $\mathfrak{q}(x) \in \mathcal{P}$, there exist a constant $b_{0}, b_{1},...,b_{n}$ such that 
 \begin{equation}\label{6.21}
 	\mathfrak{q}(x)= \sum_{k=0}^{n} b_{k}\,\,{_{\mathcal{B}} \mathcal{E}_{n}}(x;y),
 \end{equation}
where
\begin{equation}\label{6.22}
	b_{k}=\frac{1}{k!}\left<\frac{e^{{t}}+1}{2}\,e^{-y(e^{{t}}-1)}{t}^{k}\,\, \arrowvert \,\, \mathfrak{q}(x)\right>
\end{equation}
\begin{proof}
 From \eqref{6.13}, \eqref{6.16} and \eqref{6.19}, we noted that	
 \begin{equation*}
 	\left<\frac{e^{{t}}+1}{2}\,e^{-y(e^{{t}}-1)}{t}^{k}\,\, \arrowvert \,\, {_{\mathcal{B}} \mathcal{E}_{n}}(x;y)\right>= n!\delta_{n,k}\,\,\,\,\,\,\,n,k\in \mathbb{N}\cup\{0\},
 \end{equation*}
By using \eqref{6.21}, we obtain
 \begin{equation*}
 	\aligned
	\left<\frac{e^{{t}}+1}{2}\,e^{-y(e^{{t}}-1)}{t}^{k}\,\, \arrowvert \,\, \mathfrak{q}(x)\right>=&\left<\frac{e^{{t}}+1}{2}\,e^{-y(e^{{t}}-1)}{t}^{k}\,\, \arrowvert \,\,\sum_{l=0}^{n} b_{l}\,\, {_{\mathcal{B}} \mathcal{E}_{n}}(x;y)\right>\\
	=&\sum_{l=0}^{n} b_{l}\left<\frac{e^{{t}}+1}{2}\,e^{-y(e^{{t}}-1)}{t}^{k}\,\, \arrowvert \,\,\,\, {_{\mathcal{B}} \mathcal{E}_{n}}(x;y)\right>\\
	=&\sum_{l=0}^{n} b_{l}{l!}\delta_{l,k}=k!b_{k}
	\endaligned
\end{equation*}
Which is desired result \eqref{6.22}.
\end{proof}
\end{theorem}

\begin{theorem}
 If $n \in \mathbb{N}$, we have
 \begin{equation}
 	\int_{x}^{x+z}{_{\mathcal{B}} \mathcal{E}_{n}}(v;y)\,dv= \frac{e^{z{t}} - 1}{{t}}\,\,{_{\mathcal{B}} \mathcal{E}_{n}}(x;y)
 \end{equation}
\begin{proof}
By using \eqref{6.20}, we have
 \begin{equation*}
 	\aligned
	\int_{x}^{x+z}{_{\mathcal{B}} \mathcal{E}_{n}}(v;y)\,dv=& \frac{1}{n+1}\,\,\left\{{_{\mathcal{B}} \mathcal{E}_{n+1}}(x+z;y)- {_{\mathcal{B}} \mathcal{E}_{n+1}}(x;y)\right\}\\
	=& \frac{1}{n+1} \sum_{k\geq1}\binom{n+1}{k} {_{\mathcal{B}} \mathcal{E}_{n+l-k}}(x;y)\,z^{k}\\
	=& \frac{1}{n+1}\left(\sum_{k\geq0}\frac{z^{k}}{k!}\,t^{k}-1\right){_{\mathcal{B}} \mathcal{E}_{n}}(x;y)\\
	=& \frac{e^{z{t}} - 1}{{t}}\,\,{_{\mathcal{B}} \mathcal{E}_{n}}(x;y)
	\endaligned
\end{equation*}	

we get desired result.
\end{proof}
\end{theorem}

\begin{corollary}
	If $n \in \mathbb{N}\cup\{0\}$, we have 
	\begin{equation}\label{6.24}
	\int_{0}^{z}{_{\mathcal{B}} \mathcal{E}_{n}}(v;y)\,dv= \left<\frac{e^{z{t}} -1}{{t}}\arrowvert\,\, {_{\mathcal{B}} \mathcal{E}_{n}}(y)\right>
	\end{equation}
\begin{proof}
	From \eqref{6.20}, we get
	\begin{equation}\label{6.25}
		{_{\mathcal{B}} \mathcal{E}_{n}}(x;y)=\frac{{t}}{n+1}\,\,{_{\mathcal{B}} \mathcal{E}_{n+1}}(x;y)
	\end{equation}
 and by using \eqref{6.25}, we have
   \begin{equation*}
   	\aligned
   	\left<\frac{e^{z{t}} -1}{t}\arrowvert\,\, {_{\mathcal{B}} \mathcal{E}_{n}}(y)\right>=& \left<\frac{e^{z{t}} -1}{{t}}\,\arrowvert\,\, \frac{{t}}{n+1}\,\,{_{\mathcal{B}} \mathcal{E}_{n+1}}(y)\right>\\
   	=&\left<e^{z{t}} -1\,\arrowvert\,\, \frac{1}{n+1}\,\,{_{\mathcal{B}} \mathcal{E}_{n+1}}(y)\right>\\
   	=& \frac{1}{n+1}\left\{{_{\mathcal{B}} \mathcal{E}_{n+1}}(z;y)- {_{\mathcal{B}} \mathcal{E}_{n+1}}(0;y)\right\}\\
   	=& \int_{0}^{z}{_{\mathcal{B}} \mathcal{E}_{n}}(v;y)\,\,dv
   	\endaligned
   \end{equation*}
\end{proof}
\end{corollary}
For any $\mu\in \mathbb{N}\cup\{0\}$ from \eqref{3.1} the Bell based Euler polynomials of order $\mu$ are given as:
\begin{equation}\label{6.26}
	\sum_{n\geq0}{_{\mathcal{B}} \mathcal{E}_{n}^{(\mu)}}(x;y)\frac{{t}^n}{n!}=\left(\frac{2}{e^{{t}}+1}\right)^{\mu}e^{{x{t}}+y(e^{t} -1)}\,\,\,\,\,\,\, ( |{t}|<2\pi ).
\end{equation}

When ${t}$ goes to zero then ${_{\mathcal{B}} \mathcal{E}_{n}^{(\mu)}}(x;y)=1$, which means that $o\left(\left(\frac{2}{e^{{t}}+1}\right)^{\mu}e^{{x{t}}+y(e^{t} -1)}\right)=0$. Hence the generating function \eqref{6.26} of order $\mu$ is an invertible and it will be treated as Sheffer sequence.

Suppose that
\begin{equation*}
	\mathfrak{g}^{\mu}({t};y)= \frac{(e^{{t}} +1)^\mu}{2^\mu}\,e^{-y(e^{{t}} -1)}
\end{equation*}

We know that $\mathfrak{g}^{\mu}(t;y)$ is an invertible series. From \eqref{6.26} we says that   ${_{\mathcal{B}} \mathcal{E}_{n}^{(\mu)}}(x;y)$ is an appell sequence for $\mathfrak{g}^{\mu}(t;y)$. Hence from \eqref{6.17}, we have
\begin{equation*}
	{_{\mathcal{B}} \mathcal{E}_{n}^{(\mu)}}(x;y)= \frac{1}{\mathfrak{g}^{\mu}({t};y)}\,\,x^{n},
\end{equation*}

and 
\begin{equation*}
		{t}\,\,{_{\mathcal{B}} \mathcal{E}_{n}^{\mu}}(x;y)=n\,\,{_{\mathcal{B}} \mathcal{E}_{n-1}^{\mu}}(x;y)
\end{equation*}

Thus, we have
\begin{equation*}
{_{\mathcal{B}} \mathcal{E}_{n}^{\mu}}(x;y) \sim\left(\frac{(e^{{t}}+1)^{\mu}}{2^\mu}\,e^{-y(e^{{t}}-1)}, {t}\right)
\end{equation*}

Now, using above result we discuse some interesting theorem.

\begin{theorem}
If $n \geq 0$, then 
\begin{equation}\label{6.27}
	{_{\mathcal{B}} \mathcal{E}_{n}^{(\mu)}}(y)= \sum_{i_{1}+...+i_{\mu}=n}\binom{n}{i_{1},...,i_{\mu}}\,\,{_{\mathcal{B}} \mathcal{E}_{i_{\mu}}}(y)\prod_{j=1}^{\mu-1}\mathcal{E}_{i_{j}}
\end{equation}
\begin{proof}
By using \eqref{6.3} and \eqref{6.26}, we have
\begin{equation}\label{6.28}
	\left<\frac{2^{\mu}}{(e^{{t}} +1)^{\mu}}\,e^{z{t}+y(e^{{t}} -1)}\arrowvert\,\,x^{n}\right>={_{\mathcal{B}} \mathcal{E}_{n}^{(\mu)}}(z;y)=\sum_{l=0}^{n}\binom{n}{l}{_{\mathcal{B}} \mathcal{E}_{n-l}^{(\mu)}}(y)\,z^{l}
\end{equation}

and 
\begin{equation}\label{6.29}
	\aligned
	\left<\frac{2^{\mu}}{(e^{{t}} +1)^{\mu}}\,e^{z{t}+y(e^{{t}} -1)}\arrowvert\,\,x^{n}\right>=&\left<\frac{2}{e^{{t}} +1}\times...\times\frac{2}{e^{{t}} +1}e^{y(e^{{t}} -1)}\arrowvert \,x^{n}\right>\\
	=&\sum_{i_{1}+...+i_{\mu}=n}\binom{n}{i_{1},...,i_{\mu}}{_{\mathcal{B}} \mathcal{E}_{i_{\mu}}^{(\mu)}}(y)\times\mathcal{E}_{i_{1}}\times...\times\mathcal{E}_{i_{\mu-1}}.
	\endaligned
\end{equation}
From \eqref{6.28} and \eqref{6.29}, we get result \eqref{6.27}.
\end{proof}
\end{theorem}

\begin{theorem}
	If  $\,\mathfrak{q}(x)\in \mathcal{P}_{n}$, then
	\begin{equation*}
	 \mathfrak{q}(x)=\sum_{k=0}^{n}b_{k}^{\mu}\,\,{_{\mathcal{B}} \mathcal{E}_{n}^{(\mu)}}(x;y)\, \in\mathcal{P}_{n}
	\end{equation*}
where 
\begin{equation}\label{6.30}
	b_{k}^{\mu}=\frac{1}{2^{\mu}k!}\left<(e^{{t}} +1)^{\mu}e^{-y(e^{{t}} -1)}{t}^{k}\arrowvert\, \mathfrak{q}(x)\right>
\end{equation}

\begin{proof}
	Let as assume that
	\begin{equation}\label{6.31}
		\mathfrak{q}(x)=\sum_{k=0}^{n}b_{k}^{\mu}\,\,{_{\mathcal{B}} \mathcal{E}_{n}^{(\mu)}}(x;y)\, \in\mathcal{P}_{n}
	\end{equation}

with the help of \eqref{6.31}, we can write
\begin{equation}\label{6.32}
	\aligned
	\left<\left(\frac{e^{{t}}+1}{2}\right)^{\mu}\,e^{-y(e^{{t}}-1)}{t}^{k}\,\, \arrowvert \,\, \mathfrak{q}(x)\right>=&\left<\left(\frac{e^{{t}}+1}{2}\right)^{\mu}\,e^{-y(e^{{t}}-1)}{t}^{k}\,\, \arrowvert \,\,\sum_{l=0}^{n} b_{l}^{\mu}\,\, {_{\mathcal{B}} \mathcal{E}_{n}^{\mu}}(x;y)\right>\\
	=&\sum_{l=0}^{n} b_{l}^{\mu}\left<\left(\frac{e^{{t}}+1}{2}\right)^{\mu}\,e^{-y(e^{{t}}-1)}{t}^{k}\,\, \arrowvert \,\,\,\, {_{\mathcal{B}} \mathcal{E}_{n}}(x;y)\right>\\
	=&\sum_{l=0}^{n} b_{l}^{\mu}{l!}\delta_{l,k}=k!b_{k}^{\mu},
	\endaligned
\end{equation}
and we know that
\begin{equation}\label{6.33}
\left<\left(\frac{e^{{t}}+1}{2}\right)^{\mu}\,e^{-y(e^{{t}}-1)}{t}^{k}\,\, \arrowvert \,\, \mathfrak{q}(x)\right>=\frac{1}{2^{\mu}}\left<(e^{{t}} +1)^{\mu}e^{-y(e^{{t}} -1)}{t}^{k}\arrowvert\, \mathfrak{q}(x)\right>.
\end{equation}
Hence from \eqref{6.32} and \eqref{6.33}, we obtain the result \eqref{6.30}.
\end{proof}
\end{theorem}

\section{\bf Conclusions}
In this paper, we introduced the Bell based Euler polynomials of order $\alpha$ and study their various correlation, implicit sumation and derivative formula. Also, we investigate various application of Bell based Euler polynomials of order $\alpha$ by using defintion of umbral calculus. Now, it is very useful because polynomials is solution of various differential equation and also play an importnat role in multifarious area like mathematics, physics and engineering sciences.

\vspace{0.30cm}
{\bf Conflicts of Interest:}
The authors declare no conflict of interest.

\vspace{0.30cm}
{\bf Ethical approval:} This article does not contain any studies with human participants or animals performed by any of the authors.

\vspace{0.30cm}
{\bf Author Contributions:} 
All authors contributed equally to this manuscript. 

\vspace{0.30cm}
{\bf Acknowledgments:}
The authors express their deep gratitude to the anonymous referees for their critical comments and suggestions to improve this paper to its current form.

\vspace{0.30cm}
{\bf Availability of data and material:}
Data sharing is not applicable to this article as no datasets were generated or analysed during the current study.

\end{document}